\documentclass[12pt]{article}
\usepackage[cp1251]{inputenc}
\usepackage{amssymb,latexsym,amsmath}

\usepackage{amsfonts}
\usepackage[english, russian]{babel}
\usepackage{indentfirst}
\usepackage{graphicx,graphics,verbatim}
	\usepackage{cite}

\begin{document}

		УДК 517.9 
		
		MSC 26A33, 74S40

\begin{center}		
	{\bf Fractional powers of Bessel operator and its numerical calculation} \\
\end{center}	

\begin{center}
		\textbf{Durdimurod K. Durdiev}\\
	Bukhara Branch of the Institute of Mathematics at the Academy of Sciences of the Republic of Uzbekistan, Bukhara, Uzbekistan\\
	durdiev65@mail.ru\\
		\textbf{Elina L. Shishkina}\\
	Voronezh State University, Voronezh, Russia\\
	shishkina@amm.vsu.ru\\
	\textbf{Sergei M. Sitnik}\\
	Belgorod State University, Belgorod, Russia\\
	sitnik@bsu.edu.ru\\

\end{center}

\begin{center}	
Дурдимурод Каландарович Дурдиев\\
 Бухарское отделение института Математики АН РУз, Бухара, Узбекистан\\
durdiev65@mail.ru\\
Шишкина Элина Леонидовна\\
Воронежский государственный университет, Воронеж, Россия\\ 
shishkina@amm.vsu.ru\\
Ситник Сергей Михайлович\\
Белгородский государственный национальный исследовательский университет, Белгород, Россия\\
sitnik@bsu.edu.ru\\


\end{center}


{\bf Ключевые слова}: оператор Бесселя, дробные степени, преобразование Ханкеля, обобщённый сдвиг, композиционный метод, операторы преобразования


{\bf Abstract}. The article discusses the fractional powers of the Bessel operator and their numerical implementation. An extensive literature is devoted to the study of fractional powers of the Laplace operator and their applications. Such degrees are used in the construction of functional spaces, in the natural generalization of the Schr\"odinger equation in quantum theory, in the construction of the models of acoustic wave propagation in complex media (for example, biological tissues) and space-time models of anomalous (very slow or very fast) diffusion, in spectral theory etc. If we assume the radiality of the function on which the Laplace operator acts, then we receive the problem of constructing the fractional power of the Bessel operator. We propose to use a compositional method for constructing the operators mentioned earlier, which leads to constructions similar in their properties to the Riesz derivatives. The Hankel transform is considered as a basic integral transformation. On its basis, the compositional method proposed by V.V. Katrakhov and S.M. Sitnik, negative powers of the Bessel operator are constructed. The resulting operator contains the Gaussian hypergeometric function in the kernel. For further study, the generalized translation operator is considered in the article, and its properties are proved. For constructing a positive fractional power of the Bessel operator  known methods of regularization of the integral are considered. Then, a scheme for the numerical calculation of fractional powers of the Bessel operator is proposed. This scheme is based on the Taylor--Delsarte formula obtained by B.M. Levitan. Examples containing the exact and approximate values of the positive and negative powers of the Bessel operator, the absolute error, and illustrations are given. The list of references contains sources with known results on similar fractional operators, as well as applications of them.

{\bf Аннотация}. В статье рассматриваются дробные степени оператора Бесселя и их численная реализация. Обширная литература посвящена изучению дробных степеней оператора Лапласа и их приложениям. Такие степени используются при конструировании функциональных пространств,  при естественном обобщении уравнения Шрёдингера в квантовой теории, при построении моделей распространения акустических волн в сложных средах (например, биологических тканях) и пространственно-временных моделей аномальной (очень медленной или очень быстрой) диффузии, в спектральной теории и др. Если предположить радиальность функции, на которую действует оператор Лапласа, то мы придем к задаче о построении дробной степени оператора Бесселя. Мы предлагаем использовать композиционный метод построения указанных операторов, что приводит к конструкциям, близким по своим свойствам к Риссовым производным. В качестве базового интегрального преобразования рассматривается преобразование Ханкеля. На его основе композиционным методом, предложенным В. В. Катраховым и С.М. Ситником, строятся отрицательные степени оператора Бесселя. Полученный оператор содержит в ядре гипергеометрическую функцию Гаусса. Для дальнейшего изучения в статье приводится оператор обобщенного сдвига, доказываются его свойства. Известными способами достигается регуляризация интеграла при построении положительной дробной степени оператора Бесселя.  Затем предлагается схема численного вычисления дробных степеней оператора Бесселя, основанная на полученной Б.М. Левитаном формуле Тейлора--Дельсарта. Приводятся примеры, содержащие точное и приближенное значения положительной и отрицательной степеней оператора Бесселя, абсолютную погрешность и иллюстрации. В списке литературы приводятся источники, в которых содержатся известные результаты о дробных степенях операторов, рассматриваемого в статье типа, а также содержащие их приложения.

\newpage

\tableofcontents

\section{Introduction}
We will consider  fractional powers $({B}_\gamma)^{\alpha}$, $\alpha\in \mathbb{R}$ of differential Bessel operator in the form
\begin{equation}\label{Bess}
	B_\gamma= D^2+\frac{\gamma}{x}D,\qquad \gamma\geq 0, \ \ \ D:= \frac{d}{dx}.
\end{equation}
This fractional power will be constructed using Integral Transforms Composition Method (ITCM)  \cite{Fit,BookSSh,SitnikShishkinaElsevier}. 
For more detailed discussion of the 
fractional powers of \eqref{Bess} on a segment and an semi-axes we refer to 
\cite{SitSh1,SitSh2}.
Negative fractional degrees of the operator \eqref{Bess} obtained in this way are {\it one-dimensional Riesz B-potentials} studied in \cite{DAN1991,DAN1994,LyaLiz2,RAN1996,Gul1}.

The potential theory comes from mathematical physics. The most well-known areas of its application are
electrostatic and gravitational theory, probability theory,
scattering theory, biological systems and other.

First application  of  classical Riesz potentials was given by M. Riesz himself and it was a solution of  Maxwell's equations for the electromagnetic field (see \cite{Riesz1949}, p. 146, and \cite{Fremberg2}).
Maxwell's equations are fundamental equations of the classical electrodynamics and  optics. The equations completely describe all   electromagnetic phenomena in an arbitrary environment and provide a mathematical model for electric, optical and radio technologies, such as power generation, electric motors, wireless communication, lenses, radar etc. So Riesz potentials can be used for studying of realistic single particle energy levels.

An interesting fact was noticed in \cite{Caffarelli}. Namely, in this paper it was shown that  Riesz potential can be interpreted as a transmutation operator. More precisely,
the operator square root of the Laplacian was obtained from the harmonic
extension problem to the upper half space as the operator that maps the Dirichlet boundary
condition to the Neumann condition. The same result but for hyperbolic Riesz potentials was obtained in \cite{Enciso}.

\section{Construction of  fractional power of Bessel operator by ITCM}

In this section using integral transform composition method we construct the operator $\mathbf{I}^\alpha_\gamma$ 
which is the negative fractional power of Bessel operator: $ (B_\gamma)^\alpha=\mathbf{I}^\alpha_\gamma$, $\alpha>0$.
To do this, we use a suitable integral transformation it is the Hankel transform.

Let $\gamma>0$.	The one-dimentional \textbf{Hankel transform}  of  a function $f$, such that $f(x)x^\gamma {\in} L_1(0,\infty)$ is expressed as
\begin{equation}\label{Hank}
	{H}_\gamma[f](\xi)={H}_\gamma[f(x)](\xi)={f}(\xi)=\int\limits_{0}^\infty f(x)\,{j}_\frac{\gamma-1}{2}(x\xi)x^\gamma dx,
\end{equation}
where
the symbol $j_\nu$ is  used for the normalized Bessel function of the first kind (see \cite{Kipr}, p. 10, \cite{levitan}): 
\begin{equation}\label{FBess1}
	j_\nu(x) =\frac{2^\nu\Gamma(\nu+1)}{x^\nu}\,\,J_\nu(x),
\end{equation}
where  $J_\nu$ is  Bessel function of the first kind.

Let $f(x)x^\gamma {\in} L_1(0,\infty)$ and $f$ has bounded variation in a
neighborhood of a point $x$ of continuity of $f$. Then  the inversion formula
\begin{equation}\label{HankInv}
	{H}^{-1}_\gamma[\widehat{f}(\xi)](x)=f(x)=\frac{2^{1-\gamma}}{
		\Gamma^2\left(\frac{\gamma{+}1}{2}\right)}\int\limits_{0}^\infty
	{j}_\frac{\gamma-1}{2}(x\xi)\widehat{f}(\xi)\xi^\gamma\:d\xi
\end{equation}
holds.

We are looking for operator  $\mathbf{I}^\alpha_\gamma$ 
in the factorized form
\begin{equation}\label{4410}
	\mathbf{I}^\alpha_\gamma = H_{\gamma}^{-1} ( t^{-2\alpha} H_{\gamma}),
\end{equation}
where $H_\gamma$ is a Hankel transform \eqref{Hank}. It is exactly the form which is provided by ITCM.
Next we will consider  even function  $f=f(x)$ from the Schwartz space, $x\in\mathbb{R}$.

Let $x>0$. We have
$$
\left( \mathbf{I}^\alpha_\gamma f\right) (x) =	H_{\gamma}^{-1} \left[t^{-2\alpha} H_{\gamma}[f](t)\right](x)=
$$
$$
=\frac{2^{1-\gamma}}{\Gamma^2\left(\frac{\gamma+1}{2}\right)}
\int\limits_{0}^{\infty} {j}_{\frac{\gamma-1}{2}} (xt)\,
t^{\gamma-2\alpha}\,dt \int\limits_{0}^{\infty }{j}_{\frac{\gamma-1}{2}}(ty) f(y) y^\gamma dy=
$$	
$$
=\int\limits_{0}^{\infty }(xt)^{\frac{1-\gamma}{2}}J_{\frac{\gamma-1}{2} }(xt) t^{\gamma-2\alpha}\, dt
\int\limits_{0}^{\infty }(ty)^{\frac{1-\gamma}{2}}J_{\frac{\gamma-1}{2}}(ty) f(y)y^\gamma dy=
$$
$$
=x^{\frac{1-\gamma}{2}}\int\limits_{0}^{\infty }y^{\frac{\gamma+1}{2}}f(y)dy
\int\limits_{0}^{\infty } t^{1-2\alpha} J_{\frac{\gamma-1}{2} }(xt)J_{\frac{\gamma-1}{2}}(ty)dt =
$$
$$
=x^{\frac{1-\gamma}{2}}\int\limits_{0}^{x}y^{\frac{\gamma+1}{2}}f(y)dy
\int\limits_{0}^{\infty } t^{1-2\alpha} J_{\frac{\gamma-1}{2} }(xt)J_{\frac{\gamma-1}{2}}(ty)dt+
$$
$$
+x^{\frac{1-\gamma}{2}}\int\limits_{x}^{\infty }y^{\frac{\gamma+1}{2}}f(y)dy
\int\limits_{0}^{\infty } t^{1-2\alpha} J_{\frac{\gamma-1}{2} }(xt)J_{\frac{\gamma-1}{2}}(ty)dt.
$$

Using formula 2.12.31.1 from \cite{IR2} p. 209 of the form
$$
\int\limits_{0}^{\infty } t^{\beta-1}\, J_{\rho }(xt)J_{\nu }(yt)\, dt=
$$		
$$
=\left\{%
\begin{array}{ll}
$$2^{\beta-1}x^{-\nu-\beta}y^\nu
\frac{\Gamma\left(\frac{\nu+\rho+\beta}{2}\right)}{\Gamma(\nu+1)
	\Gamma\left(\frac{\rho-\nu-\beta}{2}+1\right)}
{_2F_1}\left( \frac{\nu+\rho+\beta}{2}, \frac{\nu-\rho+\beta}{2}; \nu+1; \frac{y^2}{x^2}\right) $$, & \hbox{$0<y<x$;} \\
$$2^{\beta-1}x^{\rho}y^{-\rho-\beta}
\frac{\Gamma\left(\frac{\nu+\rho+\beta}{2}\right)}{\Gamma(\rho+1)
	\Gamma\left(\frac{\nu-\rho-\beta}{2}+1\right)}
{_2F_1}\left( \frac{\nu+\rho+\beta}{2}, \frac{\beta+\rho-\nu}{2}; \rho+1; \frac{x^2}{y^2}\right)
$$, & \hbox{$0<x<y,$} \\
\end{array}%
\right.
$$
$$
x,y, {\rm Re}\,(\beta+\rho+\nu)>0;\,\, {\rm Re}\,\beta<2,
$$
where $\,_2F_1$ is hypergeometric Gauss function which  inside the circle $|z|{<}1$ determined as the sum of the hypergeometric series (see \cite{abramovic}, p. 373, formula 15.3.1)
\begin{equation}\label{FG}
	\,_2F_1(a,b;c;z)=F(a,b,c;z)=\sum\limits_{k=0}^\infty\frac{(a)_k(b)_k}{(c)_k}\frac{z^k}{k!},
\end{equation}
and for $ |z|\geq 1 $ it is obtained by analytic continuation of this series.
In (\ref{FG}) parameters $a,b,c$ and variable $z$ can be complex, and $c{\neq}0,{-}1,{-}2,{\dots}$. Multiplier
$(a)_k$
is the Pohgammer symbol is defined by
$(a)_k=a(a+1)...(a+k-1),$ $k=1,2,...,$ $(a)_0\equiv 1$.

Putting $\beta=2-2\alpha$, $\rho=\frac{\gamma-1}{2}$, $\nu=\frac{\gamma-1}{2}$ we obtain
$$
\int\limits_{0}^{\infty } t^{1-2\alpha} J_{\frac{\gamma-1}{2} }(xt)J_{\frac{\gamma-1}{2}}(ty)dt=
$$
$$
=\left\{%
\begin{array}{ll}
$$\frac{2^{1-2\alpha}y^{\frac{\gamma-1}{2}}}{x^{2-2\alpha+\frac{\gamma-1}{2}}}
\frac{\Gamma\left(\frac{\gamma+1}{2}-\alpha\right)}{\Gamma\left( \frac{\gamma+1}{2}\right)
	\Gamma\left(\alpha\right)}
{_2F_1}\left( \frac{\gamma+1}{2}-\alpha,1- \alpha; \frac{\gamma+1}{2}; \frac{y^2}{x^2}\right) $$, & \hbox{$0<y<x$;} \\
$$\frac{2^{1-2\alpha}x^{\frac{\gamma-1}{2}}}{y^{2-2\alpha+\frac{\gamma-1}{2}}}
\frac{\Gamma\left(\frac{\gamma+1}{2}-\alpha\right)}{\Gamma\left(\frac{\gamma+1}{2}\right)
	\Gamma\left(\alpha\right)}
{_2F_1}\left( \frac{\gamma+1}{2}-\alpha,1-\alpha;  \frac{\gamma+1}{2}; \frac{x^2}{y^2}\right)
$$, & \hbox{$0<x<y,$} \\
\end{array}%
\right.
$$
$$
{\rm Re}\,(\gamma+1-2\alpha)>0;
$$
and
$$
\left( \mathbf{I}^\alpha_\gamma f\right) (x) =\frac{2^{1-2\alpha}\Gamma\left(\frac{\gamma+1}{2}-\alpha\right)}{\Gamma\left( \frac{\gamma+1}{2}\right)
	\Gamma\left(\alpha\right)}\times
$$
$$
\times\left( 
\,x^{2\alpha-\gamma-1}\int\limits_{0}^{x}f(y)
{_2F_1}\left( \frac{\gamma+1}{2}-\alpha,1- \alpha; \frac{\gamma+1}{2}; \frac{y^2}{x^2}\right)y^\gamma dy+\right.
$$
$$
\left.+\,\int\limits_{x}^{\infty }f(y)
{_2F_1}\left( \frac{\gamma+1}{2}-\alpha,1-\alpha;  \frac{\gamma+1}{2}; \frac{x^2}{y^2}\right) y^{2\alpha-1}dy\right).
$$	

In \cite{Kratcer} the next formula is proved 
$$
\,_2F_1\left(a,a-b+\frac{1}{2},b+\frac{1}{2};z^2\right)=(1+z)^{-2a} \,_2F_1\left(a,b,2b;\frac{4z}{(1+z)^{2}}\right),
$$
using which we get for $a=\frac{\gamma+1}{2}-\alpha$, $b=\frac{\gamma}{2}$,
$$
x^{2\alpha-\gamma-1}y^\gamma{_2F_1}\left( \frac{\gamma+1}{2}-\alpha,1- \alpha; \frac{\gamma+1}{2}; \frac{y^2}{x^2}\right)=
$$
$$=
y^\gamma (x+y)^{2\alpha-\gamma-1}{_2F_1}\left( \frac{\gamma+1}{2}-\alpha,\frac{\gamma}{2};\gamma; \frac{4xy}{(x+y)}\right),
$$
$$
y^{2\alpha-1}{_2F_1}\left( \frac{\gamma+1}{2}-\alpha,1-\alpha;  \frac{\gamma+1}{2}; \frac{x^2}{y^2}\right)=
$$
$$=
y^\gamma (x+y)^{2\alpha-\gamma-1}{_2F_1}\left( \frac{\gamma+1}{2}-\alpha,\frac{\gamma}{2};\gamma; \frac{4xy}{(x+y)}\right).
$$
So we get a definition of negative fractional power of Bessel operator for $0{<}\alpha{<}\frac{\gamma{+}1}{2}$ 
$$
(B_\gamma)^{-\alpha}f(x)=\left( \mathbf{I}^\alpha_\gamma f\right) (x) =
$$
\begin{equation}\label{DrI}	
	=
	\frac{2^{1-2\alpha}\Gamma\left(\frac{\gamma+1}{2}-\alpha\right)}{\Gamma\left( \frac{\gamma+1}{2}\right)
		\Gamma\left(\alpha\right)}\int\limits_{0}^{\infty}f(y) (x+y)^{2\alpha-\gamma-1}{_2F_1}\left( \frac{\gamma+1}{2}-\alpha,\frac{\gamma}{2};\gamma; \frac{4xy}{(x+y)}\right)y^\gamma dy.
\end{equation}

Since the hypergeometric series \eqref{FG} converges only in the unit circle of the complex plane, it is necessary to construct an analytic continuation of the hypergeometric function beyond the boundary of this circle, to the entire complex plane. One of the ways to continue analytically is to use the Euler integral representation of the form
\begin{equation}\label{FGI}
	\,_2F_1(a,b;c;z) = { \frac{\Gamma(c)}{\Gamma(b)\Gamma(b-c)} } \int\limits_{0}^{1}
	t^{b-1} (1-t)^{c-b-1} (1-tz)^{-a} \,dt,
\end{equation}
$$
0<\textrm{Re}\,b<\textrm{Re}\,c,\qquad |\textrm{arg}(1-z)|<\pi,
$$
in which the right side is defined under the specified conditions,
ensuring the convergence of the integral.

\section{Generalized translation and fractional power of Bessel operator}

We will use  the subspace of the space of rapidly decreasing functions:
$$
S_{ev}=\left\{f\in C^\infty_{ev}:\sup_{x\in[0,\infty)}\left|x^{\alpha }D^{\beta }f(x)\right|<\infty \quad \forall \alpha ,\beta \in \mathbb {Z} _{+}\right\},
$$
where $C^\infty_{ev}$ is a class of infinetly differentiable functions on $[0,\infty)$,  such that $\frac{d^{2k+1}f}{d x_i^{2k+1}}\biggr|_{x=0}=0$ for all non-negative integer $k$.

We consider $\gamma\geq0$.  Let $L_p^{\gamma}(0,\infty)=L_p^{\gamma}$, $1{\leq}p{<}\infty$, be the space of  all measurable in $(0,\infty)$
functions, admitting even continuation on $\mathbb{R}$, such that
$$
\int\limits_{0}^\infty|f(x)|^p x^\gamma dx<\infty.
$$
For a real number $p\geq 1$, the $L_p^\gamma$--norm of $f$ is defined by
$$
||f||_{p,\gamma}=\left(\,\,\int\limits_{0}^\infty|f(x)|^p x^\gamma dx\right)^{1/p}.
$$
It is known (see \cite{Kipr}) that $L_p^\gamma$ is a Banach space. For every $1\leq p<\infty$ the Schwartz class $S_{ev}$ is dense in $L_p^{\gamma}$.

In this section we consider  the transmutation operator  called the \textit{generalized translation} of the form (see \cite{levitan}):
\begin{equation}\label{Sdvog0}
	^\gamma T^y_xf(x)=C(\gamma)\int\limits_0^\pi
	f(\sqrt{x^2+y^2-2xy\cos{\varphi}})\sin^{\gamma-1}{\varphi}d\varphi,\qquad \gamma>0,
\end{equation}
$$
C(\gamma)=\left(\int\limits_0^\pi\sin^{\gamma-1}{\varphi}d\varphi\right)^{-1}=\frac{\Gamma\left(\frac{\gamma+1}{2}\right)}{\sqrt{\pi}\,\,\Gamma\left(\frac{\gamma}{2}\right)},\qquad \gamma>0.
$$
For $\gamma=0$ generalized translation $\,^\gamma T^y_x$ is
$$
\,^0 T^y_x= T^y_xf(x)=\frac{f(x+y)+f(x-y)}{2}.
$$

It is known \cite{levitan} that
\begin{equation}\label{RavDBes}
	\,^\gamma T^y_x j_{\frac{\gamma-1}{2}}(x\xi)=j_{\frac{\gamma-1}{2}}(x\xi)\,j_{\frac{\gamma-1}{2}}(y\xi).
\end{equation}

Let consider some properties of the generalized translation.

{\bf Property 1.}
For the generalized translation operator  $\,^\gamma T^y_x$ the representation
\begin{equation}\label{Sdvig01}
	\,^\gamma T^y_xf(x)=2^{\gamma-1} C(\gamma)\int\limits_0^{1}
	f\left((x+y)\sqrt{1-\frac{4xy}{(x+y)^2}z}\,\right) z^{\frac{\gamma}{2}-1}(1-z)^{\frac{\gamma}{2}-1}dz.
\end{equation}
is valid.

\textit{Proof.}
We transform the generalized translation operator as follows. First in
\eqref{Sdvog0}
putting $\varphi=2\alpha$ we obtain
$$\,^\gamma T^y_xf(x)=2C(\gamma)\int\limits_0^{\pi/2}
f(\sqrt{x^2+y^2+2xy\cos{2\alpha}})\sin^{\gamma-1}({2\alpha})d\alpha=$$
$$=2^\gamma C(\gamma)\int\limits_0^{\pi/2}
f\left(\sqrt{x^2+y^2+2xy(\cos^2{\alpha}-\sin^2\alpha)}\right)\sin^{\gamma-1}{\alpha}\cos^{\gamma-1}{\alpha}d\alpha=$$
$$
=2^\gamma C(\gamma)\int\limits_0^{\pi/2}
f\left(\sqrt{x^2+y^2+2xy(1-2\sin^2\alpha)}\right)\sin^{\gamma-1}{\alpha}(1-\sin^2{\alpha})^{\frac{\gamma-1}{2}}d\alpha.
$$
Now let $\sin{\alpha}=t$ then ${\alpha}{=}0$ for $t{=}0$, ${\alpha}{=}\pi/2$ for $t{=}1$,
$d\alpha{=}\frac{dt}{(1-t^2)^{1/2}}$ and
$$\,^\gamma T^y_xf(x)=2^\gamma C(\gamma)\int\limits_0^{1}
f(\sqrt{x^2+y^2+2xy(1-2t^2)})t^{\gamma-1}(1-t^2)^{\frac{\gamma}{2}-1}dt=
$$
$$
=\{t^2=z\}=
$$
$$
=2^{\gamma-1} C(\gamma)\int\limits_0^{1}
f(\sqrt{x^2+y^2+2xy(1-2z)})z^{\frac{\gamma}{2}-1}(1-z)^{\frac{\gamma}{2}-1}dz=
$$
$$
=2^{\gamma-1} C(\gamma)\int\limits_0^{1}
f(\sqrt{(x+y)^2-4xyz})z^{\frac{\gamma}{2}-1}(1-z)^{\frac{\gamma}{2}-1}dz=
$$
$$
=2^{\gamma-1} C(\gamma)\int\limits_0^{1}
f\left((x+y)\sqrt{1-\frac{4xy}{(x+y)^2}z}\,\right) z^{\frac{\gamma}{2}-1}(1-z)^{\frac{\gamma}{2}-1}dz.
$$
The proof is complete.

{\bf Property 2.}
For the generalized translation operator $\,^\gamma T^y_x$ the representation
\begin{equation}\label{Sdvog011}
	(\,^\gamma T^y_xf)(x){=}\frac{2^\gamma C(\gamma)}{(4xy)^{\gamma-1}}\int\limits_{|x-y|}^{x+y}zf(z)[(z^2-(x-y)^2)((x+y)^2-z^2)]^{\frac{\gamma}{2}-1}dz
\end{equation}
is valid.

\textit{Proof.} Changing the variable $\varphi$ to the $2\alpha$ in  \eqref{Sdvog0} we obtain
$$\,^\gamma T^y_xf(x)=2C(\gamma)\int\limits_0^{\pi/2}
f(\sqrt{x^2+y^2-2xy\cos{2\alpha}})\sin^{\gamma-1}({2\alpha})d\alpha=$$
$$=2^\gamma C(\gamma)\int\limits_0^{\pi/2}
f\left(\sqrt{x^2+y^2-2xy(\cos^2{\alpha}-\sin^2\alpha)}\right)\sin^{\gamma-1}{\alpha}\cos^{\gamma-1}{\alpha}d\alpha=$$
$$
=2^\gamma C(\gamma)\int\limits_0^{\pi/2}
f\left(\sqrt{x^2+y^2-2xy(1-2\sin^2\alpha)}\right)\sin^{\gamma-1}{\alpha}(1-\sin^2{\alpha})^{\frac{\gamma-1}{2}}d\alpha.
$$
Now putting $\sin{\alpha}{=}t$ we get for ${\alpha}{=}0$, $t{=}0$, for ${\alpha}{=}\pi/2$, $t{=}1$,
$d\alpha{=}\frac{dt}{(1-t^2)^{1/2}}$ and
$$\,^\gamma T^y_xf(x)=2^\gamma C(\gamma)\int\limits_0^{1}
f(\sqrt{x^2+y^2-2xy(1-2t^2)})t^{\gamma-1}(1-t^2)^{\frac{\gamma}{2}-1}dt.
$$
Introducing the variable $z$ by the equality $\sqrt{x^2+y^2-2xy(1-2t^2)}=z$ we obtain
$$
t=\left(\frac{z^2-(x-y)^2}{4xy}\right)^{1/2},\qquad dt=\frac{zdz}{(4xy)^{1/2}(z^2-(x-y)^2)^{1/2}},
$$
$z=|x-y|$ when $t=0$,   $z=x+y$ when $t=1$ and
$$
\,^\gamma T^y_xf(x)=\frac{2^\gamma C(\gamma)}{(4xy)^{\gamma-1}}\int\limits_{|x-y|}^{x+y}zf(z)[(z^2-(x-y)^2)((x+y)^2-z^2)]^{\frac{\gamma}{2}-1}dz.
$$
The proof is complete.

{\bf Property 3.} If $f$ is a Schwartz function, $g$ is a continuous function then
\begin{equation}\label{Samos}
	\int\limits_0^{\infty}\,^\gamma T^y_xf(x)g(y)y^\gamma dy
	=\int\limits_0^{\infty}f(y)\,^\gamma T^y_xg(x)y^\gamma dy.
\end{equation}

\textit{Proof.} Applying to $\int\limits_0^{\infty}\,^\gamma T_x^yf(x)g(y)y^\gamma dy$ the representation \eqref{Sdvog011} we obtain
$$
\int\limits_0^{\infty}\,^\gamma T_x^yf(x)g(y)y^\gamma dy=
$$
$$
=(4x)^{1-\gamma}2^\gamma C(\gamma)\int\limits_0^{\infty}yg(y)dy\int\limits_{|x-y|}^{x+y}zf(z)[(z^2-(x-y)^2)((x+y)^2-z^2)]^{\frac{\gamma}{2}-1}dz=
$$
$$
=(4x)^{1-\gamma}2^\gamma C(\gamma)\biggl[\int\limits_0^{x}yg(y)dy\int\limits_{x-y}^{x+y}zf(z)[(z^2-(x-y)^2)((x+y)^2-z^2)]^{\frac{\gamma}{2}-1}dz+
$$
$$
+\int\limits_x^{\infty}yg(y)dy\int\limits_{y-x}^{x+y}zf(z)[(z^2-(x-y)^2)((x+y)^2-z^2)]^{\frac{\gamma}{2}-1}dz\biggr].
$$
Converting an expression $(z^2-(x-y)^2)((x+y)^2-z^2)$ and changing the order of integration we get
$$
\int\limits_0^{\infty} \,^\gamma T_x^yf(x)g(y)y^\gamma dy=
$$
$$
=(4x)^{1-\gamma}2^\gamma C(\gamma)\biggl[\int\limits_0^{x}zf(z)dz\int\limits_{x-z}^{x+z}yg(y)[((z+x)^2-y^2)(y^2-(z-x)^2)]^{\frac{\gamma}{2}-1}dy+
$$
$$
+\int\limits_x^{\infty}zf(z)dz\int\limits_{z-x}^{x+z}yg(y)[((z+x)^2-y^2)(y^2-(z-x)^2)]^{\frac{\gamma}{2}-1}dy\biggr]=
$$
$$
=\int\limits_0^{\infty}f(z)\,^\gamma T_x^zg(y)z^\gamma dz.
$$
The commutation is proved.

{\bf Property 4.} For $x>0$ the formula representing a generalized translation $\,^\gamma T^y_x$ of power function $x^\alpha$ is
\begin{equation}\label{TranslationPow1}
	\,^\gamma T^y_x x^\alpha =(x+y)^\alpha
	\,_2F_1\left(-\frac{\alpha}{2},\frac{\gamma}{2};\gamma;\frac{4xy}{(x+y)^2}\right),
\end{equation}
where $\,_2F_1$ is the Gaussian hypergeometric function.

\textit{Proof.}
Let first $x\neq y$.	Using formula \eqref{Sdvig01} let find $\,^\gamma T^y_x$ of  $x^\alpha$. We have
$$
\,^\gamma T^y_x
x^\alpha
=\frac{2^{\gamma-1}\Gamma\left(\frac{\gamma+1}{2}\right)}{\sqrt{\pi}\Gamma\left(\frac{\gamma}{2}\right)}(x+y)^\alpha\int\limits_0^{1}
\left(1-\frac{4xy}{(x+y)^2}z\right)^{\frac{\alpha}{2}}(1-z)^{\frac{\gamma}{2}-1}z^{\frac{\gamma}{2}-1}dz.
$$
The last integral is the Gaussian hypergeometric function \eqref{FGI} for
$z{=}\frac{4xy}{(x+y)^2}$, $a=-\frac{\alpha}{2}$, $b=\frac{\gamma}{2}$, $c=2b=\gamma$, ($c>b>0$), thus
$$
\,^\gamma T^y_x
x^\alpha=\frac{2^{\gamma-1}\Gamma\left(\frac{\gamma+1}{2}\right)\Gamma\left(\frac{\gamma}{2}\right)}{\sqrt{\pi}\Gamma(\gamma)}\,(x+y)^\alpha
\,_2F_1\left(-\frac{\alpha}{2},\frac{\gamma}{2},\gamma;\frac{4xy}{(x+y)^2}\right).
$$

Using the doubling formula for gamma function 
$$
\Gamma(2z)=\frac{2^{2z-1}}{\sqrt{\pi}}\Gamma(z)\Gamma\left(z+\frac{1}{2}\right)
$$
we obtain \eqref{TranslationPow1}.

{\bf Property 5.}	The generalized translation $\,^\gamma T^y_x$ of $e^{-x^2}$, $x>0$ is
\begin{equation}
	\label{SdExp}
	\,^\gamma T^y_x e^{-x^2}=\Gamma \left({\tfrac {\gamma+1}{2}}\right)
	\left(xy\right)^{\frac{1-\gamma}{2}}
	e^{-x^2-y^2}I_{{{\frac {\gamma-1}{2}}}}\left(2xy\right).
\end{equation}

\textit{Proof.}
Using the formula \eqref{Sdvog011} we obtain 
$$
\,^\gamma T^y_x  e^{-x^2}=\frac{2^\gamma C(\gamma)}{(4xy)^{\gamma-1}}\int\limits_{|x-y|}^{x+y}z e^{-z^2}[(z^2-(x-y)^2)((x+y)^2-z^2)]^{\frac{\gamma}{2}-1}dz.
$$
Find the integral
$$
I=\int\limits_{|x-y|}^{x+y}z e^{-z^2}[(z^2-(x-y)^2)((x+y)^2-z^2)]^{\frac{\gamma}{2}-1}dz=\{z^2=t\}=
$$
$$
=\frac{1}{2}\int\limits_{(x-y)^2}^{(x+y)^2} e^{-t}[(t-(x-y)^2)((x+y)^2-t)]^{\frac{\gamma}{2}-1}dt=\{t-(x-y)^2=w\}=
$$
$$
=\frac{1}{2}e^{-(x-y)^2}\int\limits_{0}^{4xy} e^{-w}[w(4xy-w)]^{\frac{\gamma}{2}-1}dw.
$$
Applying the formula 2.3.6.2 from \cite{IR1} of the form
\begin{equation}
	\label{IR1}
	\int\limits_0^a x^{\alpha-1}(a-x)^{\alpha-1}e^{-px}dx=\sqrt{\pi}\Gamma(\alpha)\left( \frac{a}{p}\right) ^{\alpha-1/2}e^{-ap/2} I_{\alpha-1/2}(ap/2),
\end{equation}
$$
\textrm{Re}\,\alpha>0,
$$
we get
$$
I=2^{\gamma-2}\sqrt{\pi}\Gamma\left(\frac{\gamma}{2}\right)\,e^{-x^2-y^2}(xy)^{\frac {\gamma-1}{2}}I_{\frac {\gamma-1}{2}}\left(2xy\right).
$$
Then
$$
\,^\gamma T^y_x  e^{-x^2}=\frac{2^\gamma }{(4xy)^{\gamma-1}}\frac{\Gamma\left(\frac{\gamma+1}{2}\right)}{\sqrt{\pi}\,\,\Gamma\left(\frac{\gamma}{2}\right)}\int\limits_{|x-y|}^{x+y}z e^{-z^2}[(z^2-(x-y)^2)((x+y)^2-z^2)]^{\frac{\gamma}{2}-1}dz=
$$
$$
=\frac{2^\gamma }{(4xy)^{\gamma-1}}\frac{\Gamma\left(\frac{\gamma+1}{2}\right)}{\sqrt{\pi}\,\,\Gamma\left(\frac{\gamma}{2}\right)}\,2^{\gamma-2}\sqrt{\pi}\Gamma\left(\frac{\gamma}{2}\right)\,e^{-x^2-y^2}(xy)^{\frac {\gamma-1}{2}}I_{{{\frac {\gamma-1}{2}}}}\left(2xy\right).
$$
After simplification we get \eqref{SdExp}.

Using \eqref{TranslationPow1} we can rewrite \eqref{DrI} as
$$
(B_\gamma)^{-\alpha}f(x)=\left( \mathbf{I}^\alpha_\gamma f\right) (x) =
\frac{2^{1-2\alpha}\Gamma\left(\frac{\gamma+1}{2}-\alpha\right)}{\Gamma\left( \frac{\gamma+1}{2}\right)
	\Gamma\left(\alpha\right)}\int\limits_{0}^{\infty}(\,^\gamma T_x^y x^{2\alpha-\gamma-1})f(y)y^\gamma dy.
$$

Using \eqref{Samos}	we can write
\begin{equation}
	\label{Pot}
	(B_\gamma)^{-\alpha}f(x)=\left( \mathbf{I}^\alpha_\gamma f\right) (x) =
	\frac{2^{1-2\alpha}\Gamma\left(\frac{\gamma+1}{2}-\alpha\right)}{\Gamma\left( \frac{\gamma+1}{2}\right)
		\Gamma\left(\alpha\right)}\int\limits_{0}^{\infty}(\,^\gamma T_x^y f(x)) x^{2\alpha-1} dy.
\end{equation}
It is easy to see that for Schwartz functions this integral converges for all $\alpha>0$. Expression \eqref{Pot}  we will cal {\it one-dimensional Riesz B-potential} or {\it Bessel--Riesz fractional integral}.

Let  verify that for  $\alpha=1$ and for $f=f(x)$ from the Schwartz space we get $\left( \mathbf{I}^1_\gamma B_\gamma f\right) (x)=f(x)$:
$$
(B_\gamma)^{-\alpha}B_\gamma f(x)=\left( \mathbf{I}^1_\gamma B_\gamma f\right) (x) =
\frac{\Gamma\left(\frac{\gamma-1}{2}\right)}{2\Gamma\left( \frac{\gamma+1}{2}\right)
}\int\limits_{0}^{\infty}\left( \,^\gamma T^y_x (B_\gamma)_x f(x) \right)y dy=
$$
$$
=
\frac{1}{\gamma-1}\int\limits_{0}^{\infty} \left((B_\gamma)_y \,^\gamma T^y_x f(x) \right)y dy
=\frac{1}{\gamma-1}\int\limits_{0}^{\infty} \left(\frac{d}{dy}y^\gamma \frac{d}{dy} \,^\gamma T^y_x f(x) \right)y^{1-\gamma} dy=
$$
$$
=\frac{1}{\gamma-1}\left( y^\gamma \frac{d}{dy} \,^\gamma T^y_x f(x)\right)y^{1-\gamma}\biggr|_{y=0}^{\infty} -\frac{1}{\gamma-1}\int\limits_{0}^{\infty} \left(y^\gamma \frac{d}{dy} \,^\gamma T^y_x f(x) \right)\frac{d}{dy}y^{1-\gamma} dy=
$$
$$
=\int\limits_{0}^{\infty}\frac{d}{dy} \,^\gamma T^y_x f(x)dy=\,^\gamma T^y_x f(x) \biggr|_{y=0}^\infty=f(x).
$$

For $\gamma=0$, $\alpha\neq 1,3,5,...$ and for even function $f(x)$ we get
$$
\left( \mathbf{I}^\alpha_0 f\right) (x) 
=\frac{\Gamma\left(\frac{1}{2}-\alpha\right)}{2^{2\alpha}\Gamma\left( \frac{1}{2}\right)
	\Gamma\left(\alpha\right)}\int\limits_{0}^{\infty}\left[f(x+y)-f(x-y)\right]  y^{2\alpha-1}  dy=
$$
$$
=\frac{\Gamma\left(\frac{1}{2}-\alpha\right)}{2^{2\alpha}\Gamma\left( \frac{1}{2}\right)
	\Gamma\left(\alpha\right)}\left( \int\limits_{0}^{\infty}f(x+y) y^{2\alpha-1}  dy-
\int\limits_{0}^{\infty}f(x-y) y^{2\alpha-1}  dy
\right) =
$$
$$
=\frac{\Gamma\left(\frac{1}{2}-\alpha\right)}{2^{2\alpha}\Gamma\left( \frac{1}{2}\right)
	\Gamma\left(\alpha\right)}\left( \int\limits_{-\infty}^{0}f(x-y) |y|^{2\alpha-1}  dy-
\int\limits_{0}^{\infty}f(x-y)y^{2\alpha-1}  dy
\right) =
$$
$$
=\frac{\Gamma\left(\frac{1}{2}-\alpha\right)}{2^{2\alpha}\Gamma\left( \frac{1}{2}\right)
	\Gamma\left(\alpha\right)}\int\limits_{-\infty}^{\infty}f(x-y)|y|^{2\alpha-1}  dy=
$$
$$
=\frac{\Gamma\left(\frac{1}{2}-\alpha\right)}{2^{2\alpha}\Gamma\left( \frac{1}{2}\right)
	\Gamma\left(\alpha\right)}\int\limits_{-\infty}^{\infty}f(y)|y-x|^{2\alpha-1}  dy.
$$
Since
$$
\Gamma\left(\frac{1}{2}-\alpha\right)=\frac{\pi}{\sin\left(\frac{\pi+\alpha\pi}{2}\right)\Gamma\left(\frac{1}{2}+\alpha\right) }
$$
and
$$
\Gamma(\alpha)\Gamma\left(\frac{1}{2}+\alpha\right)=\sqrt{\pi}2^{1-2\alpha}\Gamma(2\alpha) 
$$
we get the Riesz potential (see formula 12.1, p. 214 in \cite{SKM}) 
$$
\left( \mathbf{I}^\alpha_0 f\right) (x) 
=\frac{1}{2\Gamma\left( 2\alpha\right)
	\cos(\alpha\pi/2)}\int\limits_{-\infty}^{\infty}f(y)|y-x|^{2\alpha-1}  dy.
$$

Now let define Bessel--Riesz fractional derivative. To do this, 
we will use the  finite-difference method.
Due to its versatility and simplicity this method is  currently finding wider application to the problems of mathematical physics
in particular with fractional derivatives.
It is applicable both for theoretical studies of various problems, and for their approximate numerical results.
We consider generalized finite differences defined by the formula
\begin{equation}\label{necenter.o.k.r.}
	(\boxdot_{\,t}^{\,l}f)(x)=\sum\limits_{k=0}^l(-1)^kC^k_l\,^\gamma T_x^{kt}f(x).
\end{equation}

Let $l=2[\alpha]+1$, $0<\alpha$ is not entire. Namely, {\it Bessel--Riesz fractional derivative} is
\begin{equation}\label{B-g.s.i.}
	\left(\mathbf{B}_{\gamma}^\alpha
	f\right)(x)=\lim\limits_{{(L_p^\gamma)}\atop{\varepsilon\longrightarrow
			+0}}\frac{1}{d_{l,\gamma}(\alpha)}\int\limits_{\varepsilon}^\infty\frac{(\boxdot_{\,t}^{\,l}f)(x)}{t^{1+2\alpha}}\,\:dt
	=\frac{1}{d_{l,\gamma}(\alpha)}\int\limits_{0}^\infty\frac{(\boxdot_{\,t}^{\,l}f)(x)}{t^{1+2\alpha}}	\,\:dt.
\end{equation}
Here
$$
d_{l,\gamma}(\alpha)=\frac{\pi\,\Gamma\left(\gamma+1\over2\right)}{2^{2\alpha+1}\,
	\Gamma\left(\frac{1+\gamma}{2}+\alpha\right)\Gamma(\frac{1}{2}+\alpha)}
\frac{\sum\limits_{k=0}^{l}(-1)^{k+1}C_l^kk^{2\alpha}}{\sin{\alpha\pi}}.
$$
For example, for $0<\alpha<1$
$$
\left(\mathbf{B}_{\gamma}^\alpha
f\right)(x)=\frac{1}{d_{1,\gamma}(\alpha)}\int\limits_{0}^\infty\frac{f(x)-\,^\gamma T_x^{t}f(x)}{t^{1+2\alpha}}	\,\:dt.
$$

\section{Numerical scheme and examples}

Numerical methods for fractional integrals and derivatives are usually based on the Grunwald-Letnikov formulas, which are a generalization of  formulas with finite differences and Riemannian sums, or on the use of the representation of solutions by infinite series. 
In order to find Bessel--Riesz fractional integral  Gauss–Laguerre quadrature (see \cite{Davis, Rabinowitz}) can be used. Namely, the numerical integration formula of Gaussian type
$$
\int\limits_0^\infty e^{-x}f(x)dx=\sum\limits_{k=0}^n a_k f(x_k)+\frac{(n!)^2}{(2n)!}f^{(2n)}(\xi),\qquad 0<\xi<\infty.
$$
For numerical implementation, it is convenient to use recurrence formulas for
Laguerre polynomials \cite{Szego}. The first two polynomials are defined as
$$ L_{0}(x)=1,\qquad L_1(x) = 1 - x
$$
and then using the following recurrence relation for any $k\geq 1$:
$$
L_{k + 1}(x) = \frac{(2k + 1 - x)L_k(x) - k L_{k - 1}(x)}{k + 1}. 
$$

We have
$$
\int\limits_{0}^{\infty }f(y) \left( \,^\gamma T^y_x x^{2\alpha-\gamma-1}\right)y^\gamma\,dy= \sum _{i=1}^{n}w_{i}g(y_{i})+E_n(\xi),
$$
where $y_i$ is the $i$-th root of Laguerre polynomial $L_n(x)$ and the $w_i$ is given by 
$$ w_{i}={\frac {y_{i}}{\left(n+1\right)^{2}\left[L_{n+1}\left(x_{i}\right)\right]^{2}}},
$$
$$
g(y)=e^yf(y) \left( \,^\gamma T^y_x x^{2\alpha-\gamma-1}\right)y^\gamma,
$$
$$
E_n(\xi)=\frac{(n!)^2}{(2n)!}g^{(2n)}(\xi),\qquad 0<\xi<\infty.
$$
If $f(x)$ decrease exponentially as $x\rightarrow \infty$ then Gauss-Laguerre quadrature is best:
$$
E_n(\xi)<C\cdot\frac{2n+1}{4^n}.
$$
For integrands that are $O(1/x^p)$, $p > 2\alpha-1$, as $x\rightarrow \infty$, the convergence rate becomes
$$
E_n(\xi)<\frac{C}{n^{p-2\alpha+1}}.
$$

\textit{Example 1.} Let find $\mathbf{I}^\alpha_\gamma e^{-x^2}$. Using \eqref{SdExp} and formula 2.15.6.4 from \cite{IR2} we obtain
$$
\mathbf{I}^\alpha_\gamma e^{-x^2} =\frac{2^{1-2\alpha}\Gamma\left(\frac{\gamma+1}{2}-\alpha\right)}{\Gamma\left( \frac{\gamma+1}{2}\right)
	\Gamma\left(\alpha\right)}\int\limits_{0}^{\infty}\left( \,^\gamma T^y_x e^{-x^2} \right)y^{2\alpha-1} dy=
$$
$$
=\frac{2^{1-2\alpha}\Gamma\left(\frac{\gamma+1}{2}-\alpha\right)}{\Gamma\left(\alpha\right)} x^{\frac{1-\gamma}{2}} e^{-x^2}\int\limits_{0}^{\infty}e^{-y^2} I_{{{\frac {\gamma-1}{2}}}}\left(2xy\right) y^{2\alpha-1+\frac{1-\gamma}{2}} dy=
$$
$$
=\frac{\Gamma\left(\frac{\gamma+1}{2}-\alpha\right)}{2^{2\alpha}\Gamma\left( \frac{\gamma+1}{2}\right)} 
e^{-x^2}\,_1F_1\left(\alpha,\frac{\gamma+1}{2};x^2 \right). 
$$

For $n=10$, $\alpha=0.7$, $\gamma=0.5$ we get  next results.

\begin{center}
	\begin{tabular}{|c|c|c|c|} \hline 
		$x$ & Numerical calculation of  & Exact value of   & Absolute error \\ 
		& $\mathbf{I}^\alpha_\gamma e^{-x^2}$, $\alpha=0.7$, $\gamma=0.5$ &  $\mathbf{I}^\alpha_\gamma e^{-x^2}$, $\alpha=0.7$, $\gamma=0.5$ & \\
		\hline 
		0,01 & 6,020591617 & 6,020591621 & 3,84693E-09 \\ \hline 
		0,2 & 6,004767496 & 6,0047675 & 3,69646E-09 \\ \hline 
		0,39 & 5,962266113 & 5,962266116 & 3,30447E-09 \\ \hline 
		0,58 & 5,898170155 & 5,898170158 & 2,74828E-09 \\ \hline 
		0,77 & 5,81943073 & 5,819430732 & 2,1265E-09 \\ \hline 
		0,96 & 5,733357751 & 5,733357752 & 1,53078E-09 \\ \hline 
		1,15 & 5,646313464 & 5,646313465 & 1,02519E-09 \\ \hline 
		1,34 & 5,562932683 & 5,562932684 & 6,38767E-10 \\ \hline 
		1,53 & 5,485940593 & 5,485940594 & 3,7027E-10 \\ \hline 
		1,72 & 5,416426922 & 5,416426922 & 1,99684E-10 \\ \hline 
		1,91 & 5,354341757 & 5,354341757 & 1,00189E-10 \\ \hline 
		2,1 & 5,299001572 & 5,299001572 & 4,67759E-11 \\ \hline 
		2,29 & 5,249480524 & 5,249480524 & 2,03046E-11 \\ \hline 
		2,48 & 5,204851927 & 5,204851927 & 8,22364E-12 \\ \hline 
		2,67 & 5,16430296 & 5,16430296 & 3,09441E-12 \\ \hline 
		2,86 & 5,127166869 & 5,127166869 & 1,09246E-12 \\ \hline 
		3,05 & 5,0929132 & 5,0929132 & 3,58824E-13 \\ \hline 
		3,24 & 5,061123119 & 5,061123119 & 1,14575E-13 \\ \hline 
		3,43 & 5,031463853 & 5,031463853 & 4,44089E-14 \\ \hline 
		3,62 & 5,003667528 & 5,003667528 & 1,42109E-14 \\ \hline 
		3,81 & 4,977515233 & 4,977515233 & 5,32907E-15 \\ \hline 
		4 & 4,952825466 & 4,952825466 & 2,66454E-15 \\ \hline 
		4,19 & 4,929445799 & 4,929445799 & 7,99361E-15 \\ \hline 
		4,38 & 4,9072468 & 4,9072468 & 1,33227E-14 \\ \hline 
		4,57 & 4,886117511 & 4,886117511 & 3,81917E-14 \\ \hline 
	\end{tabular}

\end{center}

\begin{figure}[h!]
	\center{\includegraphics[width=0.66\linewidth]{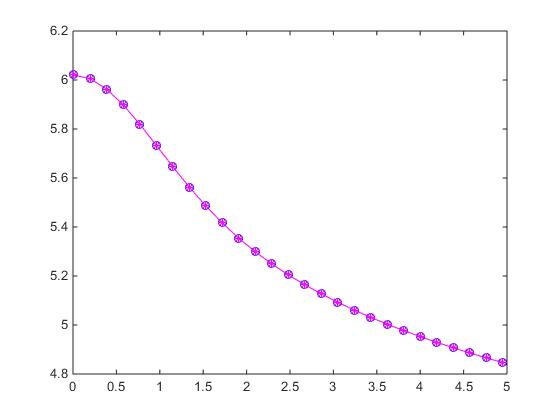} }
	\caption{$\mathbf{I}^\alpha_\gamma e^{-x^2}=\frac{\Gamma\left(\frac{\gamma+1}{2}-\alpha\right)}{2^{2\alpha}\Gamma\left( \frac{\gamma+1}{2}\right)} 
		e^{-x^2}\,_1F_1\left(\alpha,\frac{\gamma+1}{2};x^2 \right)$.}
\end{figure}

\begin{figure}[h!]
	\center{\includegraphics[width=0.66\linewidth]{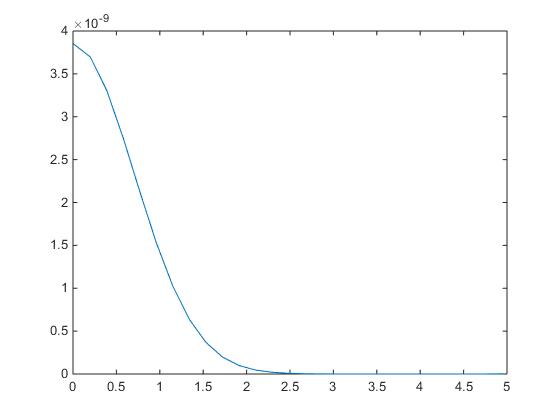} }
	\caption{Absolute error of numerical calculation of $\mathbf{I}^\alpha_\gamma e^{-x^2}$.}
\end{figure}

As for numerical scheme for  Bessel--Riesz fractional derivative we should notice that
just as ordinary shift operators can be expanded in powers of the differentiation operator, the operators $\,^\gamma T_x^y$ 
can be expanded in powers of the Bessel operator by the formula \cite{levitan}
$$
\,^\gamma T_x^y f(x)=\sum\limits_{k=0}^m\varphi_k(y) B_\gamma^k f(x)+R_m(x,y),
$$
where
$$
\varphi_k(x)=
\frac{1}{k!}\frac{\Gamma\left(\frac{\gamma+1}{2}\right)}{\Gamma\left(\frac{\gamma+1}{2}+k\right)}\left(\frac{x}{2}\right)^{2k}, 
$$
$$
R_m(x,y)=\varphi_{m+1}(y) B_\gamma^{m+1} f(\xi)|_{\xi=x+\theta y},\qquad -1<\theta<1.
$$
We have
$$
\varphi_0(x)=1,\qquad \varphi_1(x)=
\frac{\Gamma\left(\frac{\gamma+1}{2}\right)}{\Gamma\left(\frac{\gamma+1}{2}+1\right)}\left(\frac{x}{2}\right)^{2}=\frac{x^2}{2(\gamma+1)}.
$$
Then
$$
\,^\gamma T_x^{t}f(x)-f(x)=\frac{t^2}{2(\gamma+1)}B_\gamma f(\xi)|_{\xi=x+\theta t}
$$
and
for $0<\alpha<1$
$$
\left(\mathbf{B}_{\gamma}^\alpha
f\right)(x)=-\frac{1}{2(\gamma+1)d_{1,\gamma}(\alpha)}\int\limits_{0}^\infty t^{1-2\alpha}(B_\gamma f(\xi)|_{\xi=x+\theta t})	\,\:dt
$$
converges absolutely for $f\in S_{ev}$. Since $S_{ev}$ is dense in $L_p^{\gamma}$ we can also take function $f\in L_p^{\gamma}$, $1\leq p<\infty$.

\textit{Example 2.}
Let find $\mathbf{B}_{\gamma}^\alpha j_{\frac{\gamma-1}{2}}(x)$ for $\gamma=2$, $0<\alpha<1$. Using formula \eqref{RavDBes} and the fact that $j_{1/2}(x)=\frac{\sin{x}}{x}$ we get: 
$$
\left(\mathbf{B}_{\gamma}^\alpha
f\right)(x)=\frac{1}{d_{1,2}(\alpha)}\frac{\sin{x}}{x}\int\limits_{0}^\infty\frac{t-\sin{t}}{t^{2\alpha+2}}	\:dt=
$$
$$
=\Gamma(-1-2\alpha)\cos(\alpha\pi)\frac{1}{d_{1,2}(\alpha)}\frac{\sin{x}}{x}=\frac{\pi \cos(\alpha\pi)}{\Gamma(2\alpha+2)\sin(2\pi\alpha)}\frac{1}{d_{1,2}(\alpha)}\frac{\sin{x}}{x},
$$
since
$$
d_{1,2}(\alpha)=\frac{\pi\,\Gamma\left(3\over2\right)}{2^{2\alpha+1}\,
	\Gamma\left(\frac{3}{2}+\alpha\right)\Gamma(\frac{1}{2}+\alpha)}
\frac{1}{\sin{\alpha\pi}}
$$
we get
$$
\left(\mathbf{B}_{\gamma}^\alpha
f\right)(x)=\frac{2^{2\alpha}\Gamma\left(\frac{3}{2}+\alpha\right)\Gamma(\frac{1}{2}+\alpha)}{\Gamma\left(3\over2\right)\Gamma(2\alpha+2)}\frac{\sin{x}}{x}.
$$

For $n=10000$, $\alpha=0.2$, $\gamma=2$ we get  next results.

\begin{center}
	\begin{tabular}{|c|c|c|c|} \hline 
		$x$ & Numerical calculation of  & Exact value of   & Absolute error \\ 
		& $\mathbf{B}^\alpha_\gamma j_{\frac{1}{2}}(x)$, $\alpha=0.2$, $\gamma=2$ &  $\mathbf{B}^{\alpha}_\gamma j_{\frac{1}{2}}(x)$, $\alpha=0.2$, $\gamma=2$ & \\
		\hline
		0,01 & 1,397316 & 1,41372 & 0,016428 \\ \hline 
		0,3 & 1,397316 & 1,392633 & 0,016183 \\ \hline 
		0,59 & 1,397316 & 1,333139 & 0,015491 \\ \hline 
		0,88 & 1,397316 & 1,238213 & 0,014388 \\ \hline 
		1,17 & 1,397316 & 1,112569 & 0,012928 \\ \hline 
		1,46 & 1,397316 & 0,96238 & 0,011183 \\ \hline 
		1,75 & 1,397316 & 0,794917 & 0,009237 \\ \hline 
		2,04 & 1,397316 & 0,618117 & 0,007183 \\ \hline 
		2,33 & 1,397316 & 0,440132 & 0,005114 \\ \hline 
		2,62 & 1,397316 & 0,26886 & 0,003124 \\ \hline 
		2,91 & 1,397316 & 0,11151 & 0,001296 \\ \hline 
		3,2 & 1,397316 & -0,02579 & 0,0003 \\ \hline 
		3,49 & 1,397316 & -0,1383 & 0,001607 \\ \hline 
		3,78 & 1,397316 & -0,22288 & 0,00259 \\ \hline 
		4,07 & 1,397316 & -0,27812 & 0,003232 \\ \hline 
		4,36 & 1,397316 & -0,30433 & 0,003536 \\ \hline 
		4,65 & 1,397316 & -0,30344 & 0,003526 \\ \hline 
		4,94 & 1,397316 & -0,2788 & 0,00324 \\ \hline 
		5,23 & 1,397316 & -0,2349 & 0,00273 \\ \hline 
		5,52 & 1,397316 & -0,17703 & 0,002057 \\ \hline 
		5,81 & 1,397316 & -0,11089 & 0,001289 \\ \hline 
		6,1 & 1,397316 & -0,04222 & 0,000491 \\ \hline 
		6,39 & 1,397316 & 0,023587 & 0,000274 \\ \hline 
		6,68 & 1,397316 & 0,081794 & 0,00095 \\ \hline 
		6,97 & 1,397316 & 0,128612 & 0,001495 \\ \hline 
		7,26 & 1,397316 & 0,161377 & 0,001875 \\ \hline 
		7,55 & 1,397316 & 0,178666 & 0,002076 \\ \hline 
		7,84 & 1,397316 & 0,180307 & 0,002095 \\ \hline 
		8,13 & 1,397316 & 0,16731 & 0,001944 \\ \hline 
		8,42 & 1,397316 & 0,141717 & 0,001647 \\ \hline 
		8,71 & 1,397316 & 0,106388 & 0,001236 \\ \hline 
		9 & 1,397316 & 0,064737 & 0,000752 \\ \hline 
		9,29 & 1,397316 & 0,020448 & 0,000238 \\ \hline 
		9,58 & 1,397316 & -0,02281 & 0,000265 \\ \hline 
	\end{tabular}
\end{center}

\begin{figure}[h!]
	\center{\includegraphics[width=0.66\linewidth]{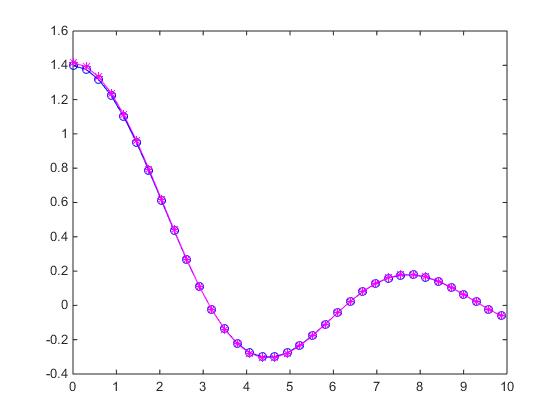} }
	\caption{$\mathbf{B}^{\alpha}_\gamma j_{\frac{1}{2}}(x)=\frac{2^{2\alpha}\Gamma\left(\frac{3}{2}+\alpha\right)\Gamma(\frac{1}{2}+\alpha)}{\Gamma\left(3\over2\right)\Gamma(2\alpha+2)}\frac{\sin{x}}{x}.$.}
\end{figure}

\begin{figure}[h!]
	\center{\includegraphics[width=0.66\linewidth]{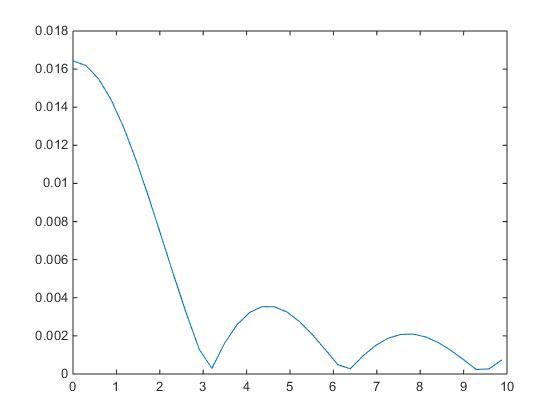} }
	\caption{Absolute error of numerical calculation of $\mathbf{B}^{\alpha}_\gamma j_{\frac{1}{2}}(x)$.}
\end{figure}

\newpage

\section{Conclusion}
On the Hankel transform basis using the integral transform compositional method proposed by V.V.Katrakhov and S.M. Sitnik, fractional powers of the Bessel operator are constructed. This article also discussed numerical methods for fractional powers of Bessel operator. Differential and integral operators presented in the article include one-dimensional Riesz B-potentials and Bessel--Riesz fractional derivative.

\end{document}